\documentclass[12pt]{article}

\usepackage{amsfonts}
\usepackage{amsmath}
\usepackage{latexsym}

\newtheorem{thm}{{\sc Theorem}}

\newtheorem{lem}{{\sc Lemma}}

\newtheorem{defn}{{\sc Definition}}

\newcommand{\qed}{\hspace*{\fill} Q.E.D.}
\newcommand{\qqed}{\hspace*{\fill} $\Box$}
 
\title{On fibre space structures of a projective 
       irreducible symplectic 
       manifold}
\author{Daisuke Matsushita 
        }
\date{E-mail Address tyler@kurims.kyoto-u.ac.jp}

\begin{document}
\maketitle
\begin{abstract}
 In this note, we investigate fibre space structures of a
 projective irreducible
 symplectic manifold. We prove that an 
 2n-dimensional projective irreducible symplectic
 manifold admits only an n-dimensional fibration 
 over a Fano variety which has only ${\mathbb Q}$-factorial
 log-terminal singularities and whose Picard number is one.
 Moreover we prove that a general fibre is an 
 abelian variety up to finite
 unramified cover, especially,
 a general fibre is an abelian surface for 4-fold.
\end{abstract}

\section{Introduction}
 We first define an {\it irreducible symplectic manifold}.
\begin{defn}
 A complex manifold $X$ is called {\it irreducible symplectic}
 if $X$ satisfies the following three conditions:
\begin{enumerate}
 \item $X$ is compact and K\"{a}hler.
 \item $X$ is simply connected.
 \item $H^{0}(X,\Omega^{2}_{X})$ is spanned by an everywhere 
       non-degenerate two-from $\omega$.
\end{enumerate}
\end{defn}
 Such a manifold can be considered as an unit of compact K\"{a}hler
 manifold $X$ with $c_1 (X) = 0$ due to the following 
 Bogomolov decomposition theorem.
\begin{thm}[Bogomolov decomposition theorem \cite{bogomolov}]
 A compact K\"{a}hler manifold $X$ with $c_{1}(X) = 0$ admits
 a finite unramified covering of $\tilde{X}$ which is
 isomorphic to a product 
 $T \times X_1 \times \cdots \times X_r \times A$
 where $T$ is a complex torus, $X_i$ are irreducible symplectic
 manifolds and $A$ is a projective manifold with 
 $h^{0}(A,\Omega^{p}) = 0$, $0 < p < \dim A$. 
\end{thm}
 In dimension 2, $K3$ surfaces are the only irreducible symplectic
 manifolds, and irreducible symplectic manifolds are considered
 as higher-dimensional analogies of $K3$ surfaces.
 In this note, we investigate {\it fibre space structures} of
 a projective irreducible symplectic manifolds.
\begin{defn}
 For an algebraic variety $X$, 
 a fibre space structure of $X$
 is a proper surjective morphism $f : X \to S$ which satisfies 
 the following two conditions:
\begin{enumerate}
 \item $X$ and $S$ are normal varieties such that  $0 < \dim S < \dim X$
 \item A general fibre of $f$ is  connected.
\end{enumerate}
\end{defn}

 Some of $K3$ surface $S$ has a fibre space structure 
 $f : S \to {\mathbb P}^{1}$ 
 whose general fibre is an elliptic curve.
 In higher dimensional analogy, we obtain the following results.
\begin{thm}
 Let $f : X \to B$ be a fibre space structure of a
 projective irreducible
 symplectic $2n$-fold $X$ with
 projective base $B$. Then
 a general fibre $F$ of $f$ and $B$ satisfy the
 following three conditions:
\begin{enumerate}
 \item $F$ is an abelian variety 
       up to finite unramified cover and $K_F \sim {\cal O}_{F}$.
 \item $B$ is $n$-dimensional and has only ${\mathbb Q}$-factorial 
       log-terminal singularities 
 \item $-K_B$ is ample and Picard number $\rho (B)$ is one.
\end{enumerate}
 Especially, if $X$ is $4$-dimensional, 
 a general fibre of $f$ is an abelian
 surface.
\end{thm}
{\sc Example. \quad}
 Let $S$ be a $K3$ surface with an elliptic fibration 
 $g: S \to {\mathbb P}^{1}$ and $S^{[n]}$
 a $n$-pointed Hilbert scheme of $S$. It is known that
 $S^{[n]}$ is an irreducible symplectic $2n$-fold and
 there exists a birational morphism 
 $\pi : S^{[n]} \to S^{(n)}$ where $S^{(n)}$ is the 
 symmetric $n$-product of $S$ (cf. \cite{beauville}). 
 We can consider 
 $n$-dimensional abelian fibration
 $ g^{(n)} : S^{(n)} \to {\mathbb P}^{n}$ for the symmetric $n$-product
 of $S^{(n)}$. 
 Then the composition morphism 
 $g^{(n)} \circ \pi : S^{[n]} \to {\mathbb P}^{n}$
 gives an example of a fibre space structure of an
 irreducible symplectic manifold.

\vspace{5mm}
\noindent
{\sc Remark. \quad} 
 Markushevich obtained some result of theorem 2 in 
 \cite[Theorem 1, Proposition 1]{mark1}
 under the assumption $\dim X = 4$ and 
 $f : X \to B$ is the moment map.
 In general, a fibre space structure of 
 an irreducible symplectic manifold is 
 not a moment map.
 Markushevich constructs in \cite[Remark 4.2]{mark2} counterexample.

\vspace{5mm}
\noindent
{\sc Acknowledgment. \quad} 
 The author express his thanks to Professors Y.~Miyaoka, S.~Mori 
 and N.~Nakayama for their advice and encouragement. He also
 thanks to Prof. D.~Huybrechts \cite{hyubrechts} for his nice survey articles 
 of irreducible symplectic manifolds.
 
\section{Proof of Theorems}
 First we introduce the following theorem due to
 Fujiki \cite{fujiki} and Beauville \cite{beauville}.
\begin{thm}[\cite{fujiki} Theorem 4.7, Lemma 4.11, Remark 4.12
            \cite{beauville} Th\`{e}or\'{e}me 5]
 \label{quadtatic}
 Let $X$ be an irreducible symplectic $2n$-fold. 
 Then there
 exists a nondegenerate quadratic form $q_{X}$ of
 signature $(3,b_{2}(X) - 3)$ on 
 $H^{2}(X , {\mathbb Z})$ which satisfies
\begin{eqnarray*}
  \alpha^{2n} &=& a_0 q_{X}(\alpha , \alpha)^{n} \\
  c_{2i}(X)\alpha^{2n - 2i} &=& a_i q_{X}(\alpha , \alpha )^{n-i} \quad 
    (i \ge 1),  
\end{eqnarray*}
 where $\alpha \in H^{2}(X, {\mathbb Z})$ and
 $a_i$'s are constants depending on $X$.
\end{thm} 
 We shall prove theorem 2 in five steps.
\begin{enumerate}
 \item $\dim B = n$ and $B$ has only log-terminal singularities;
 \item A general fibre $F$ of $f$ is an abelian variety
       up to unramified finite cover and $K_F \sim {\cal O}_{F}$;
 \item $\rho (B) = 1$;
 \item $B$ is ${\mathbb Q}$-factorial;
 \item $-K_B$ is ample.
\end{enumerate}
\noindent
{\sc Step 1. \quad} $\dim B = n$ and $B$ has only log-terminal
                    singularities.

\begin{lem}\label{pseado}
 Let $X$ be 
 an irreducible symplectic projective $2n$-fold and
 $E$ be a divisor on $X$ such that 
 $E^{2n} = 0$.
 Then,
\begin{enumerate}
 \item If $E.A^{2n-1} = 0$ 
       for some ample divisor $A$, $E \equiv 0$.
 \item If  $E.A^{2n-1} > 0$ for an ample divisor $A$  
       on $X$, then
$$
\left\{
\begin{array}{ccc}
   E^{m}A^{2n-m}  & = 0 & (m > n) \\
                  & > 0 & (m \le n)
\end{array}
\right.
$$
\end{enumerate} 
\end{lem} 
{\sc Proof of lemma. \quad}
 Let $V := \{ E \in H^{2}(X , {\mathbb Z}) | E.A^{2n-1} = 0 \}$.
 By \cite[Lemma 4.13]{fujiki}, $q_{X}$ is negative definite
 on $W$ where $V = H^{2,0} \oplus H^{0,2} \oplus W$.
 Thus, if $E.A^{2n-1} = 0$ and $E^{2n} = 0$, $E \equiv 0$. 
 Next we prove (2).
 From Theorem \ref{quadtatic}, for every integer $t$,
\begin{equation}\label{key}
 (tE + A)^{2n} = a_0 (q_{X}(tE+A , tE+A))^n .
\end{equation}
 Because $E^{2n} = a_0 (q_{X}(E,E))^n = 0$, 
$$
 q_{X}(tE+A, tE+A) = 2tq_{X}(E,A) + q_{X}(A,A).
$$
 Thus the right hand side of the equation (\ref{key})
 has order at most $n$. Comparing the both hand side of 
 the equation (\ref{key}),
 we can obtain $E^{m}.A^{2n-m} = 0$ for $m > n$. 
 If $E.A^{2n-1} > 0$, comparing the first order term of $t$ of both hand 
 of the equation (\ref{key}) we can obtain $q_{X}(E,A) > 0$.
 Because coefficients of other terms of 
 left hand side of (\ref{key}) can be written $q_{X}(E,A)$ 
 and $q_{X}(A,A)$, we can obtain
 $E^{m}.A^{2n-m} > 0$ for $0 < m \le n$.
\qqed

\noindent
 Let $H$ be a very ample divisor on $B$. Then $f^{*}H$ is a 
 nef divisor such that $(f^{*}H)^{2n} = 0$, $(f^{*}H).A^{2n-1} > 0$
 for an ample divisor $A$ on $X$.
 Thus $\dim B = n$.
 From \cite[Theorem 2]{nakayama}, $B$ has only log-terminal
 singularities.

\vspace{5mm}
\noindent 
{\sc Step 2. \quad} A general fibre $F$ of $f$ is an abelian variety
       up to unramified finite cover and $K_F \sim {\cal O}_{F}$. 

\vspace{5mm}
\noindent
 By adjunction, $K_F \sim 0$. Moreover
$$
 c_2 (F) = c_2 (X)(f^{*}H)^{2n - 2} = a_1 (q_{X}(f^{*}H,f^{*}H))^{n-1} = 0,
$$
 by Theorem \ref{quadtatic}.
 Thus $F$ has an \'{e}tale cover $\tilde{F} \to F$ such
 that $\tilde{F}$ is an Abelian variety
 by \cite{yau}.

\vspace{5mm}
\noindent
{\sc Step 3. \quad} $\rho (B) = 1$.
\begin{lem}
 Let $E$ be a divisor of $X$ such that $E^{2n} = 0$ and
 $E^{n}.(f^{*}H)^{n} = 0$. Then $E \sim_{{\mathbb Q}} \lambda f^{*}H$
 for some rational number $\lambda$.
\end{lem} 
{\sc Proof of lemma. \quad} 
 Considering the following equation
\begin{eqnarray*}
 (E - \lambda f^{*}H)^{2n} &=& 
  a_0 q_{X}(E - \lambda f^{*}H , E - \lambda f^{*}H )^n \\
    &=& a_0 (2\lambda q_{X}(E,f^{*}H))^{n},
\end{eqnarray*}
 we can obtain $q_{X}(E,f^{*}H) = cE^{n}.(f^{*}H)^{n} = 0$
 where $c$ is a constant. Thus
 $(E - \lambda f^{*}H)^{2n} = 0$.  
 Because $f^{*}H . A^{2n-1} > 0$ for every ample divisor $A$ on $X$,
 we can take a rational number $\lambda$ such that
 $(E - \lambda f^{*}H ).A^{2n-1} = 0 $ 
 Then $E -\lambda f^{*}H \equiv 0$
 by lemma \ref{pseado}.
\qqed

\noindent 
 Let $D$ be a Cartier divisor on $B$. Then $(f^{*}D)^{2n} = 0$
 and $(f^{*}D)^{n}.(f^{*}H)^{n} = 0$, thus $E \sim_{{\mathbb Q}} 
 \lambda H$ and $\rho (B) = 1$.

\vspace{5mm}
\noindent
{\sc Step 4. \quad} $B$ is ${\mathbb Q}$-factorial.
 
\vspace{5mm}
\noindent
 Let $D$ be an irreducible and reduced Weil divisor on $B$
 and $D_i$, $(1 \le i \le k)$ divisors on $X$ whose supports 
 are contained in $f^{-1}(D)$. 
 We construct  a divisor
 $\tilde{D} := \sum \lambda_i D_i$, $(\tilde{D} \not\equiv 0)$
 such that 
 $\tilde{D}^{2n} = 0$. Let $A$ be a very ample divisor on $X$,
 $S := A^{n-1}.(f^{*}H)^{n-1}$ and $C := H^{n-1}$. Then
 there exists a surjective morphism $f' : S \to C$. 
 If we choose $H$ and $A$ general, we may assume that
 $S$ and $C$ are smooth and
 $C \cap D$ are contained smooth locus of $B$.
 Because $D$ is a Cartier divisor in a neighborhood 
 of $C \cap D$, we can define $f^{'*}D$ in a neighborhood $U$
 of $S$. We can express $f^{'*}D = \sum \lambda_i D_i$ in $U$ 
 and let $\tilde{D} := \sum \lambda_i D_i$.
 Note that if $\lambda_i > 0$, $f(D_i) = D$ because
 we choose $C$ generally.
 Compareing the $n$th order term of $t$ 
 of the both hand side of the following equotion 
\begin{eqnarray*}
 (\tilde{D} + tf^{*}H )^{2n} &=& 
    a_0 q_{X}(\tilde{D}+ t f^{*}H, \tilde{D} + t f^{*}H )\\
    &=& a_0 (q_{X}(\tilde{D},\tilde{D}) + 
            2t q_{X}(\tilde{D},f^{*}H) )^{n},
\end{eqnarray*}
 we can see that 
 $ \tilde{D}^{n}.(f^{*}H)^{n} =
 cq_{X}(\tilde{D},f^{*}H)$.
 Since $f(\tilde{D}) = D$,
 $\tilde{D}^{n}.(f^{*}H)^{n} = 0$ 
 and $q_{X}(\tilde{D},f^{*}H) = 0$. 
 Considering the following equation
\begin{eqnarray*}
 (s\tilde{D} + tA + f^{*}H)^{2n} &=& 
 a_0 q_{X}( s\tilde{D} + tA + f^{*}H , s\tilde{D} + tA + f^{*}H)^n \\
 &=& a_0 ( s^2 q_{X}(\tilde{D},\tilde{D}) 
          + t^2 q_{X}(A,A) + 2stq_{X}(\tilde{D},A) \\
 & &      + 2t q_{X}(A,f^{*}H))^{n},
\end{eqnarray*}
 we can obtain 
 $q_{X}(\tilde{D},\tilde{D})q_{X}(A,f^{*}H) 
 = c \tilde{D}^2 .A^{n-1}. (f^{*}H)^{n-1}$
 where $c$ is a constant. 
 Since $\tilde{D}.A^{n-1}.(f^{*}H)^{n-1}$ is a fibre
 of $f'$, $\tilde{D}^2 .A^{n-1}.(f^{*}H)^{n-1} = 0$.
 Thus $a_0 q_{X}(\tilde{D},\tilde{D}) 
 = \tilde{D}^{2n} = 0$.
 Considering $\tilde{D}^{n}(f^{*}H)^{n} = 0$,
 we can obtain $\tilde{D} \sim_{\mathbb Q} 
 \lambda f^{*}H$ by Lemma \ref{pseado}, and
 $D \sim_{{\mathbb Q}} \lambda H$ because $f(\tilde{D}) = D$.
 Therefore $B$ is ${\mathbb Q}$-factorial.

\vspace{5mm}
\noindent
{\sc Step 5. \quad} $-K_B$ is ample.

\vspace{5mm}
 From Step 3,4, we can write $-K_B \sim_{{\mathbb Q}} tH$.
 It is enough to prove $t > 0$.
 Because $K_X \sim {\cal O}_{X}$ and a general fibre of
 $f : X \to B$ is a minimal model, $\kappa (B) \le 0$ 
 by \cite[Theorem 1.1]{kawamata} and $t \ge 0$. 
 Assume that $t = 0$. If $K_B \not\sim {\cal O}_{B}$,
 we can consider the following diagram:
$$
 \begin{array}{ccccc}
   & X & \to & B & \\
   \alpha & \uparrow & & \uparrow & \beta \\
   & \tilde{X} & \to & \tilde{B} & ,
 \end{array}
$$
 where $\beta$ is an unramified finite cover and 
 $K_{\tilde{B}} \sim {\cal O}_{\tilde{B}}$. 
 Because
 $\pi_{1}(X) = \{1 \}$, $\tilde{X}$ is the direct 
 sum of $X$.
 Thus there exists a morphism from $X$ to $\tilde{B}$
 and we may assume that $K_B \sim {\cal O}_{B}$.
 Then there exists a holomorphic $n$-form $\omega'$ on $X$ coming from $B$.
 However, if $n$ is odd, it is a contradiction
 because there exist no holomorphic $(2k-1)$-form on $X$.
 If $n$ is even, it is also a contradiction because
 $\omega'$ dose not generated by $\omega \in H^{0}(X , \Omega^{2})$.
 Thus $t > 0$ and we completed the proof of Theorem 2.
\qed

\begin{flushleft}
 Resarch Institute of Mathematical Science, \\
 Kyoto University. \\ 
 KITASHIRAKAWA, OIWAKE-CHO, \\
 KYOTO, 606-01, JAPAN.

\end{flushleft}

\begin{thebibliography}{9}
 \bibitem{beauville} A.~Beauville, 
                     {\it Vari\'{e}t\'{e}s k\"{a}hlerinennes dont
                          la premi\`{e}re classes de Chern est nulle}\/,
                     J. Diff. Geom., {\bf 18} (1983), 755--782.  
 \bibitem{bogomolov} F.A.~Bogomolov, 
                  {\it On the decomposition theorem of K\"{a}hler
                       manifolds with trivial canonical class}\/,
                  Math. USSR-Sb., {\bf 22} (1974), 580--583.
 \bibitem{fujiki} A.~Fujiki, {\it On the de Rham Cohomology Group
                      of a Compact K\"{a}hler Symplectic Manifold}\/,
               in {\it Algebraic geometry, 
              Adv. Stud. Pure Math.}\/, vol {\bf 10}, (ed. T.~Oda),
              Kinokuniya and North-Holland (1987), 105--165.
 \bibitem{hyubrechts} D.~Huybrechts,
                 {\it Compact HyperK\"{a}hler Manifolds: Basic
                  Results}\/, alg-geom9705025.
 \bibitem{kawamata} Y.~Kawamata,
                 {\it Minimal models and the Kodaira dimension 
                      of algebraic fibre spaces}\/,
                  J. Reine. Angew. Math., {\bf 363} (1985), 1--46.
 \bibitem{mark1} D.G.~Markushevich,  
                {\it Completely integrable projective symplectic
                 4-dimensional varieties}\/,
                Izv. Math., {\bf 59}:1 (1995), 159--187. 
 \bibitem{mark2} {\leavevmode\hbox to 3em{\hrulefill}\,}
                {\it Integrable symplectic structures on compact complex
                     manifolds}\/,
                 Math. USSR-Sb., {\bf 59}:2 (1988), 459--469.
 \bibitem{nakayama} N.~Nakayama, 
                {\it The Singularity of the Canonical Model
                     of Compact K\"{a}hler Manifolds}\/,
                 Math. Ann., {\bf 280} (1988), 509--512.
 \bibitem{yau} S.T.~Yau,
                {\it Calabi's conjecture and some new results in
                     algebraic geometry}\/,
               Proc. Nat. Acad. Sc. U.S.A., {\bf 74} (1977), 
               1978--1979.
\end{thebibliography}
\end{document}